\documentclass[12pt]{amsart}
\usepackage{amsmath}
\usepackage{amscd}
\usepackage{graphics}
\usepackage{latexsym}

\oddsidemargin .25in \theoremstyle{plain}

\newtheorem{Thm}{Theorem}
\newtheorem{Lem}[Thm]{Lemma}

\newtheorem{Prop}[Thm]{Proposition}
\newtheorem{Def}[Thm]{Definition}

\newcommand{\tb}{\mathop{\rm tb}\nolimits}

\newcommand{\bfr}{{\mathbb{R}}}
\newcommand{\OB}{\mathfrak{{ob}}}
\newcommand{\Li}{{\mathbb{L}}}

\def\v{\vskip.12in}

\begin{document}

\title[]{Embedding fillings of contact 3-manifolds}
\author{Burak Ozbagci}
\address{Department of Mathematics \\ Ko\c{c} University \\ Istanbul, Turkey}
\email{bozbagci@ku.edu.tr}

\begin{abstract}
In this survey article we describe different ways of embedding
fillings of contact 3-manifolds into closed symplectic 4-manifolds.
\end{abstract}

\maketitle \setcounter{section}{-1}
\section{Introduction}

\noindent One of the most exciting news regarding the topology of
3-manifolds in 2004 was the solution of the ``Property P''
conjecture by Kronheimer and Mrowka \cite{km-04}. Namely, they
proved that no surgery on a knot in $S^3$ can produce a
counter-example to the Poincar\'e conjecture. The last ingredient in
their proof was supplied by a recent theorem of Eliashberg
\cite{el-04}: \emph{Any weak filling of a contact 3-manifold can be
embedded symplectically into a closed symplectic 4-manifold}. This particular way
of embedding a weak filling into a closed symplectic
4-manifold was also used by Ozsv\'ath and Szab\'o \cite{ozsz-04} to show that their
(appropriately twisted) contact Heegaard Floer invariant of a fillable contact
structure does not vanish.

\v
\noindent
In order to prove his theorem
Eliashberg attaches a symplectic 2-handle
along the binding of an open book compatible
with the given weakly fillable contact structure such that
the other end of the cobordism given by this symplectic 2-handle attachment
symplectically  fibres over $S^1$.
Then he fills in this symplectic fibration by a symplectic
Lefschetz fibration over $D^2$ to obtain a symplectic
embedding of a weak filling into
a closed symplectic 4-manifold. Note that the method of
construction in \cite{el-04} takes its roots from the one
considered in \cite{ao-02}.

\v
\noindent
Eliashberg's theorem was obtained independently by Etnyre
\cite{et-04} using different methods. The first step in Etnyre's construction is to embed
a weak filling into a weak filling of an integral homology sphere.
Note that, from the \emph{surgery} point of view, this step also
fairly easily follows from Stipsicz's results in \cite{st-03}.
Then one can modify the symplectic form near
the boundary so that it becomes a strong filling (\cite{el-91}, \cite{ohon-99}).
This is just a homological argument. Now the problem is reduced to finding an
embedding of a strong filling.
The strategy here is to find a concave filling to cap off the convex boundary
of this strong filling from the ``other side''.
One can reduce the problem (cf. \cite{eh-02}) to the
existence of a symplectic
embedding of a Stein filling into a closed symplectic 4-manifold,
which was already provided by Lisca and  Matic \cite{lm-97}.
Alternatively, one can proceed with constructing
an explicit concave filling (cf. \cite{gay-02}) obtained by
a careful investigation of the monodromies of the open books compatible
with different types of symplectic and contact surgeries.

\v
\noindent
The purpose of this survey
article is to describe and compare embeddings due to Eliashberg
and Etnyre and discuss some previous work on the subject. We note that
there are now many ways of embedding a weak filling symplectically into a
closed symplectic 4-manifold. In Section~\ref{last}
we construct an embedding which is obtained by combining the various ideas
developed in the article. We would
like to point out that these embeddings are constructed by making
use of a recent theory developed by Giroux \cite{gi-02} which establishes a
(one-to-one) correspondence between open book decompositions of
3-manifolds and contact structures.

\section{Open book decompositions and contact structures}

\noindent We will assume throughout this paper that a
contact structure $\xi=\ker \alpha$ is coorientable
(i.e., $\alpha$ is a global 1-form) and positive
(i.e., $\alpha \wedge d\alpha >0 $ )
unless otherwise stated.
In the following we describe the compatibility of an
open book decomposition with a given contact structure on a
3-manifold.

\v

\noindent Suppose that for a link $L$ in a $3$-manifold $Y$ the complement $Y \setminus L$
fibers as $\pi \colon Y \setminus L \to S^1$ such that the fibers are interiors of
Seifert surfaces of $L$. Then $(L, \pi)$ is an \emph{open book
decomposition} (or just an \emph{open book}) of $Y$. For each $t \in S^1$,
the Seifert surface $F={\overline {\pi^{-1}(t)}}$ is called
a \emph{page}, while $L$
the \emph{binding} of the open book. The monodromy
of the fibration $\pi$ is called the \emph{monodromy} of the open
book decomposition.

\v

\noindent Any locally trivial bundle with fiber $F$ (a compact oriented surface)
over an oriented circle
is canonically isomorphic to the fibration $$\frac{I\times F}{(1,x)\sim
\big (0,h(x)\big )} \to \frac{I}{\partial I} \approx S^1$$ for some
self-diffeomorphism $h$ of $F$. In fact, the map $h$ is determined
by the fibration up to isotopy and conjugation by an orientation
preserving self-diffeomorphism of $F$. The isotopy class
represented by the map $h$ is called the monodromy of the
fibration.

\v
\noindent The \emph{mapping class group} $\Gamma_F$ of $F$
is defined as the quotient of the group of orientation
preserving
self-diffeomorphisms of $F$ fixing $\partial F$
pointwise modulo isotopies fixing $\partial F$
pointwise. Given a compact oriented surface $F$ with
nonempty boundary and $h \in \Gamma_F$, then we can form the mapping torus
$F(h)=I\times F/(1,x)\sim \big
(0, h(x)\big )$. Since $h$ is the identity on $\partial F$, the
boundary $\partial F(h)$ can be
canonically identified with $r$ copies of $T^2 = S^1 \times S^1$,
where the first $S^1$ factor is identified with $I/\partial I$ and
the second one is identified with a component of $\partial F$. Hence
$F(h)$ can be completed to a closed
$3$-manifold $Y$ equipped with an open book decomposition by
gluing in $r$ copies of $D^2\times S^1$ to $F(h)$ so that
$\partial D^2$ is identified with $S^1 = I/\partial I$ and the
$S^1$ factor in $D^2 \times S^1$ is identified with a boundary
component of $\partial F$. In
conclusion, an element $h\in \Gamma_F$ determines a $3$-manifold
together with an open book decomposition on it.

{\Thm [Alexander \cite{al-23}] Every closed and oriented 3-manifold admits an
open book decomposition.}

\v

\noindent The contact
condition $\alpha \wedge d \alpha >0$
can be strengthened in the presence of an open
book decomposition on $Y$ by requiring that $\alpha>0 $ on the binding and
$d\alpha>0 $ on the pages.

{\Def \label{compatible} An open book decomposition of a
$3$-manifold $Y$ and a contact structure $\xi$ on $Y$ are called
\emph{compatible} if $\xi$ can be represented by a contact form
$\alpha$ such that the binding is a transverse link, $d \alpha$ is
a symplectic form on every page and the orientation of the transverse
binding induced by $\alpha$ agrees with the boundary orientation
of the pages.}

{\Thm [Giroux \cite{gi-02}] \label{giroux}
Every contact 3-manifold admits a compatible open
book (with a connected binding).}

\v

\noindent We refer the reader to \cite{et} and \cite{os-04} for more
on the correspondence between open books and contact structures.

\section{Different types of fillings of contact 3-manifolds}

\noindent In this section we give definitions of different types
of symplectic fillings of contact 3-manifolds. A symplectic 4-manifold $(X, \omega)$
will be assumed to be oriented by $\omega \wedge \omega$.

\subsection{Weak filling} A contact $3$-manifold $(Y,\xi )$ is said to be
 \emph{weakly fillable} if there is a compact symplectic
4-manifold $(W, \omega )$ such that $\partial W= Y$ as oriented
manifolds and ${\omega|}_{\xi } > 0$. In this case we say that
$(W, \omega )$ is a \emph{weak filling} of $(Y,\xi )$.

\subsection{Strong filling} A contact $3$-manifold $(Y,\xi )$ is said to be
 \emph{strongly fillable} if there is a compact symplectic
4-manifold $(W, \omega )$ such that $\partial W= Y$ as oriented
manifolds, $\omega $ is exact near the boundary and its primitive
$\alpha$ (i.e., a $1$-form with $d\alpha =\omega $) can be chosen
in such a way that $\ker ({\alpha |}_{\partial W})=\xi$. In this
case we say that $(W, \omega )$ is a \emph{strong filling} of
$(Y,\xi )$. Clearly a strong filling is a weak filling by definition.

\v
\noindent Suppose
that $(W, \omega)$ is a compact symplectic 4-manifold with
nonempty boundary $\partial W =Y$ and there exists a Liouville
vector field $v$ ( i.e., $\mathcal{L}_v \omega =\omega $) defined in a neighborhood of  and transverse to
$Y$. Then $v$ induces a contact structure $\xi= \ker \alpha$ on
$Y$ where $\alpha=\iota_v\omega|_Y$ is a contact 1-form. If $v$ points out of $W$ along $Y$
then we say that $(W, \omega)$ is a \emph{convex filling} of $(Y, \xi)$, and $(Y, \xi)$ is said
to be the \emph{convex boundary}
of $(W, \omega)$. It is easy to see that the notion of a convex  filling is
the same as the notion of a strong filling. If $v$
points into $W$ along $Y$, on the other hand,  then we say
that $(W, \omega)$ is a \emph{concave filling} of $(Y,
\xi)$ and $(Y, \xi)$ is said
to be the \emph{concave boundary}
of $(W, \omega)$. Here notice
that if $v$ points out of $W$ then $\xi$ is a
positive contact structure on $Y$, while if $v$ points into $W$ then
$\xi$ is a positive contact structure on $-Y$.

\v
\noindent If a boundary component $Y$ of a compact symplectic
4-manifold $W$ which has a disconnected boundary satisfies the definition of
convexity (concavity, resp.) above then we say that $Y$ is a convex (concave, resp.)
boundary component of $W$.

\subsection{Stein filling} A compact 4-manifold $W$ with nonempty
boundary $\partial W =Y$
is called a \emph{Stein domain} if there
is a Stein surface $X$ with plurisubharmonic function $\varphi\colon X\to [0, \infty ) $
such that $W= \varphi^{-1} \big([0 , t]\big)$ for some regular value $t$. So a compact
manifold with boundary (and a complex structure $J$ on its interior) is a Stein domain if it
admits a proper plurisubharmonic function $\varphi$ which is constant on the boundary.
Then the complex line distribution induced by $J$ is a contact
structure $\xi$ on $Y$. In this case we say that the contact
$3$-manifold $(Y, \xi )$ is \emph{Stein fillable} and $(W,
J)$ is a called a \emph{Stein filling} of $(Y,\xi )$. It is easy to verify that a Stein filling
is a strong filling. In fact, $ d J^* (d \varphi) $ induces a K\"{a}hler structure on $(W, J)$. More
generally, a cobordism $W$ (with boundary $-Y_1\cup Y_2$) is a \emph{Stein cobordism} if
$W$ is a complex cobordism with a plurisubharmonic function $\varphi\colon W \to \bfr $ such that
$\varphi^{-1}(t_i)=Y_i$, for $t_1<t_2$.

\section{Embedding a Stein filling} \label{steinembed}

\noindent The first result in the literature about embedding a filling of
a contact 3-manifold into
a closed symplectic 4-manifold was obtained by Lisca and Matic.
Recall that a Stein filling (i.e., a Stein domain) admits a K\"{a}hler
form $ d J^* (d \varphi) $ which is an exact symplectic form, where
$\varphi$ is the plurisubharmonic function defining the Stein filling.

{\Thm [Lisca-Matic \cite{lm-97}] \label{liscamatic} A Stein
filling admits a K\"{a}hler embedding into a (minimal) compact
K\"{a}hler surface $X$ (of general type), such that the pull-back of
the K\"{a}hler form on $X$ is the exact symplectic form on the Stein filling. }

\v
\noindent
Apparently what motivated Lisca and Matic to construct such an
embedding was their search for a method to distinguish tight contact
structures. Using
Seiberg-Witten theory coupled with their embedding result,
Lisca and Matic were able to show that for
any positive integer $n$, there exists a homology 3-sphere with at
least $n$ homotopic but non-isomorphic tight contact structures.
Lisca and Matic use analytical tools in the construction of their
embedding and the starting point of their embedding is given by a
holomorphic embedding of a Stein domain into an affine algebraic
manifold with trivial normal bundle (cf. \cite{dls-94}). Roughly
speaking, the idea here is to approximate analytical maps by
algebraic ones, namely by polynomials.

\v
\noindent
A very different approach to embed a Stein filling \emph{smoothly} into
a closed symplectic 4-manifold was presented in \cite{ao-02}. The
construction in \cite{ao-02} is topologically more explicit than
the method of Lisca and Matic although the result is weaker since only
the smoothness of the embedding is clear from the presentation.

\v
\noindent
The simple construction in \cite{ao-02} is based on a theorem of Loi
and Piergallini (\cite{lp-01}, cf. also \cite{ao-01}) which says that
every Stein domain
admits a PALF, i.e., a positive allowable Lefschetz fibration over
$D^2$. It
is easy to see that the boundary of a PALF has a canonical open
book decomposition and we can assume that the binding of this open book is connected.
To embed a Stein filling (which has a
PALF structure) into a closed symplectic 4-manifold, we first attach a
2-handle to the binding of the open book in the boundary of the
PALF to get a Lefschetz fibration over $D^2$ (with closed fibers).
Then we extend this fibration  to a Lefschetz fibration over $S^2$.
The resulting
4-manifold is known to be symplectic by a result of Gompf (\cite{gs-99}).
This construction gives a smooth embedding of a Stein filling into a closed
symplectic 4-manifold.

\v
\noindent
It turned out
that the 2-handle above can be attached symplectically along the binding
(cf. \cite{el-04}) and hence the embedding in \cite{ao-02} is in
fact a symplectic embedding rather than just a smooth embedding.
We will discuss the details of the symplectic 2-handle attachment
along the binding of an open book in Section~\ref{weakembed}.

\section{Embedding a strong filling}\label{strongembed}
\noindent In  \cite{eh-02}, Etnyre
and Honda proved that every contact 3-manifold has (infinitely
many distinct) concave fillings. Their proof was based on
the embedding result of Lisca and Matic we discussed in the previous
section. In \cite{gay-02}, Gay proved the same
existence result (independent of the Lisca-Matic embedding)
by presenting a method to explicitly construct,
handle by handle, a concave filling of a given contact
3-manifold. A symplectic  embedding of a strong filling of a contact 3-manifold
into a closed symplectic 4-manifold trivially  follows from
Proposition~\ref{conca}.

\v

{\Prop [Etnyre-Honda \cite{eh-02}, Gay \cite{gay-02}] \label{conca}
Any contact 3-manifold admits a concave filling.}

\begin{proof} We will describe a proof (cf. \cite{os-04})
which is very similar to the one given in \cite{eh-02}.
The difference here is that we rather do not
translate contact $(\pm 1)$-surgeries along Legendrian knots into
the monodromy language of open books.

\v

\noindent Given an arbitrary
contact 3-manifold $(Y, \xi)$.
Let $\alpha$ be a contact 1-form for $\xi$. Consider a compact piece
$(W = Y \times I, \omega = d(t \alpha))$  of
the symplectization of $(Y, \xi)$.
It is easy to see that $Y \times \{1\}$ is a
convex boundary component of $(W, \omega)$
while $Y \times \{0\}$ is a concave boundary component.
Our strategy will be to  cap off the convex boundary of $(W, \omega)$
obtaining a concave filling of $Y \times \{0\}$.

\v

\noindent In \cite{dg-04}, Ding and Geiges  proved that every (closed)
contact 3-manifold $(Y, \xi)$ can be given by a contact $(\pm
1)$-surgery on a Legendrian link $\Li$ in the standard contact
$S^3$. Here the surgery coefficients are measured with respect to
the contact framing. Let $\Li^{\pm} \subset \Li$ denote the
set of the components of the link $\Li$ with
$(\pm 1)$-surgery coefficients, respectively.  Let ${\mathbb{K}}$
denote the Legendrian link we get by considering Legendrian
push-offs of the components of $\Li^+$.

{\Prop [Weinstein \cite{we-91}] \label{strong}
Let $(W, \omega)$ be a compact symplectic 4-manifold with a convex boundary component
$(Y, \xi)$. A 2-handle can be attached symplectically to $(W, \omega)$
along a Legendrian knot $L \subset (Y, \xi)$ in such
a way that the symplectic structure extends to the
2-handle and the new symplectic $4$-manifold $(\widetilde{W}, \widetilde{\omega})$
has a convex boundary component $(\widetilde{Y}, \widetilde{\xi})$, where
$(\widetilde{Y}, \widetilde{\xi})$ is given
by contact $(-1)$-surgery (i.e., Legendrian surgery)
along $L\subset (Y, \xi )$.}

\v

\noindent Thus  when we attach symplectic
2-handles to $(W, \omega )$ along the knots of
${\mathbb{K}} \subset Y=Y\times \{1\}$
we get a symplectic 4-manifold $(W', \omega ')$ with a convex boundary
component $(Y', \xi ')$ by Proposition~\ref{strong}. We observe that
the contact manifold $(Y', \xi ')$ can be given by a  Legendrian surgery along
$\Li^-$, since a combination of a contact $(+1)$-surgery on a Legendrian
knot in $\Li^+$ and a contact $(-1)$-surgery on its push-off in ${\mathbb{K}}$
cancels out (cf. \cite{dg-04}).
We note that the cancellation of these contact $(\pm 1)$-surgeries
just corresponds to the cancellation of
a right-handed Dehn
twist along a curve with a left-handed Dehn twist along a
curve parallel to it in the
monodromy of an open book in the proof of Etnyre and Honda \cite{eh-02}.

\v
\noindent Consequently $(Y', \xi ')$ is Stein fillable by a result of Eliashberg \cite{el-90}
since it is obtained from the standard contact $S^3$ via Legendrian surgeries only.
Consider a Stein filling
$(W'',J)$ of $(Y', \xi ')$ and embed this filling into a closed
symplectic $4$-manifold $(Z, \omega_Z)$ using Theorem~\ref{liscamatic}.
Then since a Stein filling is a convex filling by definition,
$(Z \setminus \mbox{int\;} W'')$ will be a concave filling of $(Y', \xi ')$. Hence we conclude that
$$(W', \omega') \; \bigcup_{Y'} \; (Z \setminus \mbox{int\;} W'' , \omega_Z)$$ is a
concave filling of $(Y, \xi )$, which is illustrated in Figure~\ref{str}.
Here we use Lemma~\ref{gluing} to glue these symplectic
4-manifolds symplectically.

\begin{figure}[ht]
  \begin{center}
     \includegraphics{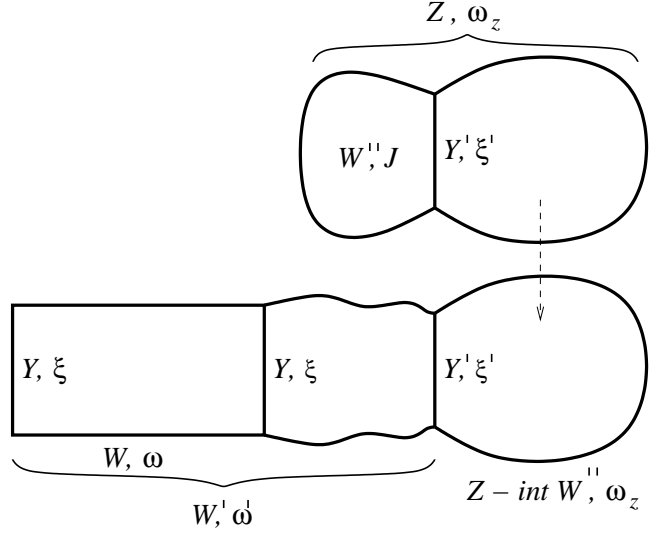}
   \caption{A concave filling of $(Y, \xi)$} \label{str}
    \end{center}
  \end{figure}

\end{proof}
\v
\noindent Next we will discuss another proof of Theorem~\ref{conca}
 given in \cite{et} which is not based
on the embedding of Lisca and Matic. This method of proof is essentially
due to Gay (\cite{gay-02}) except
for a slight short-cut at the end. We first collect below a few results that
we will need. We denote by $t_{\beta}$ a right-handed Dehn twist about a curve $\beta$
on a surface $F$.

\begin{figure}[ht]
  \begin{center}
     \includegraphics{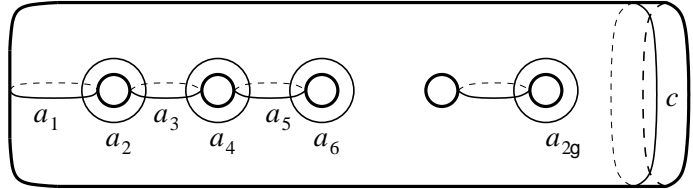}
   \caption{Genus $g$ surface $F$ with boundary} \label{surface}
    \end{center}
  \end{figure}

\v
{\Lem [Wajnryb \cite {wa-99}] \label{rel}
The relation $t_c = ( t_{a_1} \circ  t_{a_2}
\circ \cdots \circ  t_{a_{2g-1}} \circ t_{a_{2g} })^{4g+2} $
holds in the mapping class group $\Gamma_F$, where $a_i $'s are
the curves on a genus $g$ surface $F$ with one boundary component
depicted in Figure~\ref{surface} and $c$
is a curve parallel to $\partial F$. }

\v
{\Lem \label{factor} Any element $\phi$ of the mapping class group of a surface $F$ with
one boundary component can be expressed as $\phi = t_c^m \circ t_{\gamma_1}^{-1}\circ \cdots \circ
t_{\gamma_n}^{-1}$ for some $m \in \mathbb{Z}$ and some non-separating curves $\gamma_i \subset F$, where
$c$ is a curve parallel to $\partial F$.}

\begin{proof} This lemma follows easily from Lemma~\ref{rel}. It is well-known that
the mapping class group of a surface is generated by Dehn twists.
Consider an arbitrary factorization of $\phi \in \Gamma_F$ into Dehn twists. Then we can express $t_{a_1}$
---and hence any right-handed non-separating Dehn twist---by a product of non-separating
left-handed Dehn twists and $t_c$ by Lemma~\ref{rel}. Similarly applying
Lemma~\ref{rel} to subsurfaces of $F$, and using what we proved so far,  we can express
a separating Dehn twist as a product of non-separating left-handed Dehn twists and arbitrary
twists along $c$.
\end{proof}

\v
{\Lem [Honda \cite{ho-00}] \label{lrp}
A non-separating curve $\gamma$ on a convex surface in a
contact 3-manifold can be made Legendrian by isotoping the surface
through convex surfaces such that the contact framing of $\gamma$
agrees with its surface framing.}

\v
{\Lem  [Gay \cite{gay-02}] \label{compa} Given a Legendrian knot $L$ on
a page of an open book $\OB_{\xi}$
compatible with $(Y, \xi)$. Let $h \in \Gamma_F$ denote the monodromy of $\OB_{\xi}$.
Then a contact $(-1)$-surgery on $L$
induces a contact structure $\xi'$ compatible with the open
book $\OB_{\xi'}$ whose monodromy is given by $h' = h \circ t_L \in \Gamma_F$.}

\v
{\Lem [Etnyre \cite{et-98}] \label{gluing} If $(W_1, \omega_1)$ is a convex filling of
$(Y, \xi)$ and $(W_2, \omega_2)$ is a concave filling of
$(Y,\xi)$ then $X= W_1 \cup W_2$ has a symplectic structure
$\omega$ such that $\omega|_{W_1} = \omega_1$ and $ \omega|_{W_2 \setminus
\partial W_2} = k\omega_2$ where $k$ is a positive constant.}

\v
\noindent Now we are ready to describe a second proof of Theorem~\ref{conca}.
Consider the compact piece $(W, \omega)$ of the symplectization
of  $(Y, \xi)$ as in the proof above.
Let $\OB_\xi$ be an open book compatible with $\xi$ with a connected binding.
Let $\phi$ be the monodromy of this open book. Now use Lemma~\ref{factor} to write
$ \phi = t_c^m \circ  t_{\gamma_1}^{-1}\circ \cdots \circ
t_{\gamma_n}^{-1}.$ We can assume that the curve $\gamma_n$
lies on a convex page of $\OB_\xi$. Then Legendrian Realize
$\gamma_n$ (cf. Lemma~\ref{lrp})
on this convex page and perform contact $(-1)$-surgery on the Legendrian
$\gamma_n$. Performing contact $(-1)$-surgery along $\gamma_n$ induces a contact structure
which has a compatible open book whose monodromy is given by $ \phi\circ t_{\gamma_n} = t_c^m \circ
t_{\gamma_1}^{-1} \circ \cdots \circ t_{\gamma_{n-1}}^{-1}$ (cf. Lemma~\ref{compa}).

\v
\noindent We repeat this process for
all the curves $\gamma_i$ (for $i=n-1, \cdots, 1$) to obtain an open book $\OB_{\xi'}$
with monodromy $t_c^m$ which is compatible with the
resulting contact 3-manifold $(Y', \xi')$ induced by the
successive contact $(-1)$-surgeries. Moreover we can assume that $m > 1$, otherwise
we can just perform some more contact $(-1)$-surgeries along $a_i$'s depicted in Figure~\ref{surface}
(cf. Lemma~\ref{rel}).

\v
\noindent On the other hand, by Proposition~\ref{strong}, a contact $(-1)$-surgery
along a Legendrian knot $L$ in a convex boundary component of a symplectic 4-manifold
can be obtained by a symplectic 2-handle attachment along $L$. Hence
there exists a symplectic 4-manifold $(W', \xi')$ with a convex boundary component
$(Y', \xi')$ which is
obtained from $(W, \omega)$ by
attaching symplectic 2-handles along $\gamma_i$'s in the
convex end of $(W, \omega)$. Next we will prove that we can actually assume that $m=1$.
We can assume that $m$ is odd by attaching extra symplectic
2-handles along $a_i$'s if necessary.

\v
\noindent We note that Proposition~\ref{strong}
is also true
for attaching symplectic 1-handles. Namely,
one can attach a symplectic 1-handle to a sympectic 4-manifold with a convex boundary
component in such a way that the symplectic structure extends over the
1-handle and the new symplectic $4$-manifold has a convex boundary component.
At the contact 3-manifold level, the induced surgery on the convex boundary component
corresponds to taking a connected sum with
a copy of standard contact $S^1 \times S^2$. One can also view
this surgery at the level of compatible open books: Attaching a (4-dimensional)
symplectic 1-handle to a convex boundary component (along the binding of a compatible open book)
corresponds to attaching a (2-dimensional) 1-handle to the page of that open book.
Note that we extend the old monodromy by
identity over the new 1-handles.

\begin{figure}[ht]
  \begin{center}
     \includegraphics{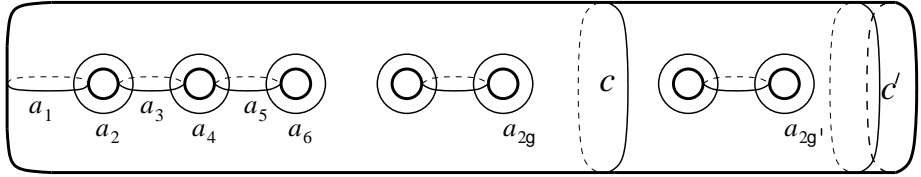}
   \caption{Genus $g'$ surface $F'$ with boundary} \label{surfbig}
    \end{center}
  \end{figure}

\v
\noindent
Now we attach symplectic 1-handles to $(W', \xi')$ so that
the resulting compatible open book on the boundary has a
page $F'$ of genus $g' = mg + \frac{m-1}{2}$ with one boundary component
so that $ m(4g+2)= 4g'+2$. Let $c'$ be a
curve parallel to $\partial F'$ as shown in Figure~\ref{surfbig}.

\v
\noindent Then by Lemma~\ref{rel} we have $$t_{c'} = ( t_{a_1} \circ t_{a_2}
\circ \cdots \circ  t_{a_{2g'-1}} \circ t_{a_{2g'} })^{4g'+2}$$
$$= ( t_{a_1} \circ  t_{a_2}
\circ \cdots \circ t_{a_{2g}} \circ t_{a_{2g+1}} \circ \cdots \circ
 t_{a_{2g'-1}} \circ t_{a_{2g'} })^{m(4g+2)} \in \Gamma_{F^{\prime}}.$$

\v
\noindent
To simplify the notation we will denote the result
of attaching symplectic 1-handles to $(W', \omega')$ again as $(W', \omega')$.
Now attach more symplectic 2-handles to $(W', \omega')$ along the curves
$a_{2g+1}, a_{2g+2}, \cdots a_{2g'}$ sufficiently many times so that
the resulting convex boundary has a compatible open book with monodromy $t_{c'}$.
Here note that we are inserting (rather
than appending as in Lemma~\ref{compa}) some right-handed Dehn twists,
but nevertheless Lemma~\ref{compa} holds true in this case (cf. \cite{gay-02}).
We will still denote the resulting symplectic 4-manifold
by $(W', \omega')$, to simplify the notation.

\v
\noindent
Summarizing the above discussion, by attaching symplectic
1- and 2-handles to $(W, \omega)$ we end up with a symplectic
4-manifold $(W', \omega')$ with a convex boundary component
$(Y', \xi')$ whose compatible open book $\OB_{\xi'}$
has the following description: The page $F'$ is a genus $g'$ surface
with one boundary component and the monodromy is a single
right-handed Dehn twist along a curve $c'$ parallel to $\partial F'$.
Let $\widehat{F}$ denote the surface obtained by
capping off the surface $F'$ by gluing a
2-disk along $\partial F'$.
\v
\noindent

{\Lem  [\cite{ao-02}] \label{circle}The 3-manifold $Y'$ is a circle bundle over
the surface $\widehat{F}$ with Euler number $-1$.}

\begin{proof}
Recall the relation $$
{(t_{a_1} t_{a_2} \cdots t_{a_{2g'}})}^{4g'+2} =t_{c'}
$$
in the mapping class group $\Gamma_{F'}$. It induces a relation $$
{(t_{a_1} t_{a_2} \cdots t_{a_{2g'}})}^{4g'+2} =1.
$$
in the mapping class group $\Gamma_{\widehat{F}}$.
This later  relation induces a Lefschetz fibration
$f : X \to S^2$ admitting a section of square $-1$.
Consider a
neighborhood $U$ of a regular fiber union this section.
We observe that $\partial U = -Y $. This is because
$X \setminus \mbox{int\;} U $ is a PALF with monodromy $$( t_{a_1} \circ  t_{a_2}
\circ \cdots \circ  t_{a_{2g'-1}} \circ t_{a_{2g'} })^{4g'+2} = t_{c'}. $$
Moreover $U$ is obtained by plumbing
a $ D^2 \times \widehat{F}$ and a disk bundle over $S^2$ with Euler number $-1$.
In Figure~\ref{diskbun1} we illustrated the handlebody diagram of the 4-manifold $U$.

\begin{figure}[ht]
  \begin{center}
     \includegraphics{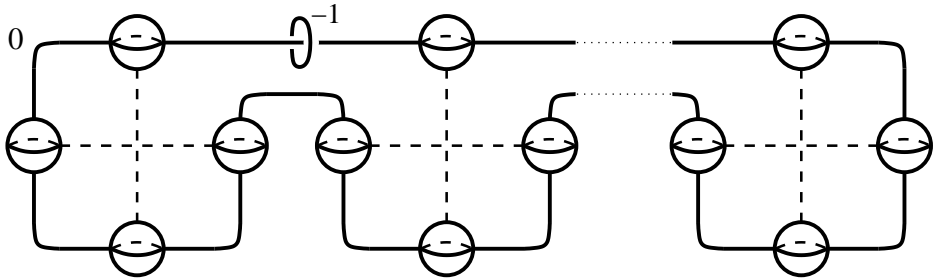}
   \caption{Plumbing a $ D^2 \times \widehat{F}$ and a $D^2$-bundle over $S^2$
   with Euler number $-1$} \label{diskbun1}
    \end{center}
  \end{figure}

\v
\noindent
We can blow down the $-1$ sphere to get a disk bundle over
$\widehat{F}$ with Euler number $+1$ (cf. Figure~\ref{diskbun2}).
Blowing down a $-1$ sphere changes the 4-manifold but the boundary 3-manifold
remains the same (up to diffeomorphism). Note that
the boundary of a disk bundle over $\widehat{F}$ with
Euler number $+1$ is circle bundle over
$\widehat{F}$ with Euler number $+1$.
Our claim follows by reversing the orientations, since when we change
the orientation of a circle bundle over
$\widehat{F}$ with Euler number $+1$, we get a circle bundle over
$\widehat{F}$ with Euler number $-1$.

\begin{figure}[ht]
  \begin{center}
     \includegraphics{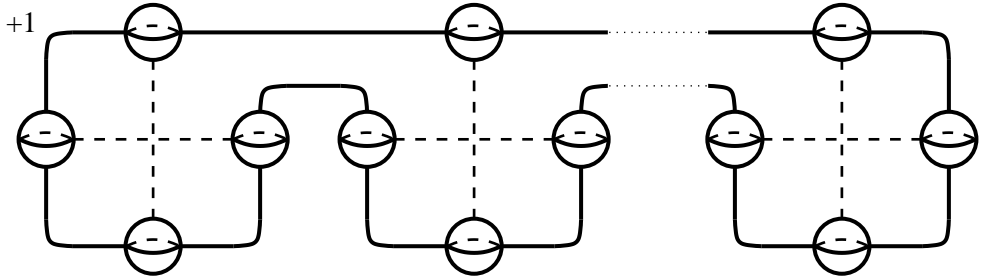}
   \caption{$D^2$-bundle over $\widehat{F}$
   with Euler number $+1$} \label{diskbun2}
    \end{center}
  \end{figure}

\end{proof}

\v
\noindent Now consider the disk bundle $M$ over $\widehat{F}$
with Euler number $1$. Then
$M$ admits a natural symplectic structure $\omega_M$
so that $(M, \omega_M)$  has a concave boundary $(Y', \xi')$.
Thus $$(W', \omega') \; \bigcup_{Y'} \; (M, \omega_M) $$
is a concave filling of $(Y, \xi)$ by Lemma~\ref{gluing}.
This finishes the proof of Proposition~\ref{conca}.

\v
\noindent Finally we would like to point out how Gay's proof in
\cite{gay-02} differs from the proof in \cite{et}.
Consider the open book $\OB_{\xi'}$ which is compatible with $(Y',\xi')$ as above.
Then Gay explains how to attach a symplectic 2-handle along the binding of
$\OB_{\xi'}$ with framing $+1$ relative to the page framing of the binding
so that the resulting contact 3-manifold $(Y'', \xi'')$ is also a \emph{concave} boundary component
of the
symplectic cobordism given by the 2-handle attachment. Denote the resulting
symplectic 4-manifold obatined by attaching this 2-handle to $(W', \omega')$
by $(W'', \omega'')$. Moreover the monodromy of the open book compatible with $(Y'', \xi'')$
is given by the identity map.  This implies that $(Y'', \xi'')$ is
contactomorphic to the standard tight contact $(S^1 \times S^2, \xi_{st})$.  Note that there is a standard
convex filling $(S^1 \times D^3, \omega_{st})$ of $(S^1 \times S^2, \xi_{st})$. Hence
$$  (W''  , \omega'' ) \; \bigcup \; (S^1 \times D^3, \omega_{st}) $$
is a concave filling of
$(Y, \xi)$ by Lemma~\ref{gluing}.

{\Prop  [Etnyre-Honda \cite{eh-02}, Gay \cite{gay-02} ] \label{stembed}

If\/ $(W, \omega )$ is a strong filling of\/ $(Y, \xi )$ then $W$
can be symplectically embedded into a closed symplectic
$4$-manifold.}

\begin{proof}
Suppose that $(W, \omega )$ is a strong filling of $(Y, \xi )$.
Consider a concave filling $(W_1, \omega_1)$ of $(Y, \xi)$.
Then we can glue (cf.
Lemma~\ref{gluing}) the symplectic manifolds $(W, \omega)$ and
$(W_1, \omega_1)$ along their common boundary $(Y, \xi )$ to get
a closed symplectic 4-manifold including $(W, \omega )$ as a
symplectic subdomain.

\end{proof}

\section{Embedding a weak filling}\label{weakembed}

\noindent In this section we will give the most general embedding result
that will cover the cases in Sections~\ref{steinembed} and
~\ref{strongembed}.

{\Thm [Eliashberg \cite{el-04}, Etnyre \cite{et-04}] \label{wembed} If\/ $(W, \omega
)$ is a weak filling of\/ $(Y, \xi )$ then $W$ can be
symplectically embedded into a closed symplectic $4$-manifold.}

\subsection{Eliashberg's construction:} Let $(W, \omega)$
be a weak filling of a contact 3-manifold $(Y, \xi )$ and let $\OB_\xi$ be an
open book decomposition of $Y$ (with a connected binding $B$)
compatible with the contact structure $\xi$ (cf. Theorem~\ref{giroux}).
Attach a symplectic
2-handle along $B$ to $Y \times I$ to obtain a symplectic
cobordism with boundary $-Y \cup Y'$ such that $Y'$ fibers over $S^1$ with
symplectic fibers. Then fill in $Y' \to S^1$ by a symplectic
Lefschetz fibration over $D^2$ to complete $W \cup H$ into a closed
symplectic 4-manifold.

\v

\noindent Eliashberg's idea above is to reduce the
question of embedding a weak filling  to a question of embedding a
symplectic surface fibration over the circle. Notice that the
binding $B$ of $\OB_\xi$ is transverse to $\xi$,  so the crucial point
of Eliashberg's construction is the way that he attaches a
symplectic 2-handle along the transverse binding $B$. Eliashberg's
construction is ``topologically" equivalent to the construction
that was given in \cite{ao-02} to embed a Stein filling
\emph{smoothly} into a closed symplectic 4-manifold. Now we
proceed with the details of Eliashberg's construction.

\v

\noindent We first describe the symplectic 2-handle $H$ to be attached along
the transverse binding $B$. We identify $\mathbb{C}^2 (z_1, z_2)$
with $\mathbb{R}^4 (x_1, y_1, x_2, y_2)$ as usual: $z_1 = x_1
+iy_1 $ and $z_2 = x_2 +iy_2$. Let $(r_i, \varphi_i)$ denote the
polar coordinates in the $z_i$-plane for $i=1,2$. Then the
standard symplectic 2-form $\omega_0$ on $\mathbb{R}^4$ is given
by
$$\omega_0 = dx_1 \wedge dy_1 + dx_2 \wedge dy_2 = r_1dr_1 \wedge
d \varphi_1  + r_2 dr_2 \wedge d \varphi_2.$$

\v

\begin{figure}[ht]
  \begin{center}
     \includegraphics{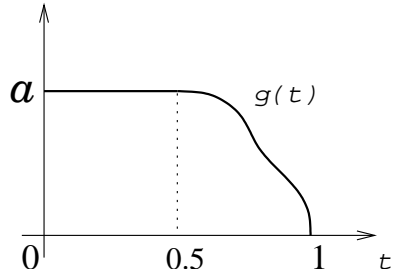}
   \caption{The graph of the smooth function $g$}
    \label{g}
    \end{center}
  \end{figure}

\begin{figure}[ht]
  \begin{center}
     \includegraphics{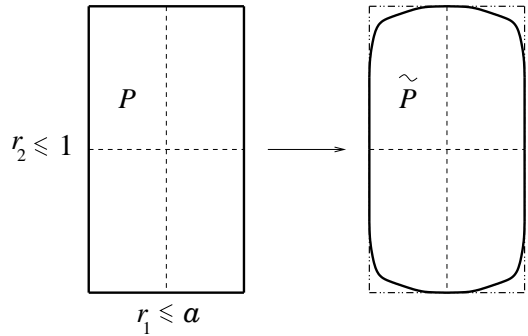}
   \caption{Smoothing the corners of the polydics}
   \label{corners}
    \end{center}
  \end{figure}

\noindent Let $a$ be a positive real number and let $P= \{ r_1 \leq a , r_2
\leq 1 \} \subset \mathbb{C}^2 $ be a polydisc. Now we define a
domain $\widetilde{P} = \{r_1 \leq g(r_2) : r_2  \in [0,1]\}
\subset P$ for some non-increasing smooth function $g(t): [0,1] \to
[0,a] $ as in Figure~\ref{g}, where $g([0, 0.5]) = a $ and $g'(t) < 0 $
for $t \in (0.5, 1)$. We will determine the real
number $a$ and the particular form of the function $g(t)$ near
$t=1$ along the proof. Here we can view $\widetilde{P}$ obtained from the polydics $P$ by
smoothing its corners as in Figure~\ref{corners}.

\noindent

\begin{figure}[ht]
  \begin{center}
     \includegraphics{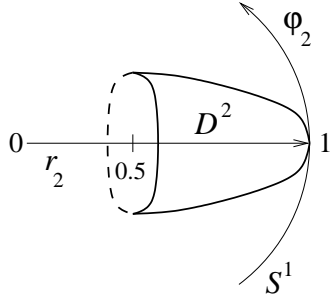}
   \caption{$\Gamma \simeq S^1 \times D^2$}
   \label{solidtorus}
    \end{center}
  \end{figure}

\v
\noindent We define $\Gamma = \{r_1=g(r_2) : r_2 \in [0.5 , 1] \}$
as part of $\partial \widetilde{P}$.
(There is a typo here in \cite{el-04}, $r_2 \in [0.5 , 1]$ not $r_1$.) We observe that $\Gamma$ is
diffeomorphic to $S^1 \times D^2$: As $r_2$ increases from
$0.5$ to $1$ in the $z_2$-plane (with polar coordinates
$(r_2, \varphi_2)$) the boundary of the disks in the $z_1$-plane
will shrink smoothly (according to the function $g$) to a point
tracing out a disk for each fixed $\varphi_2$. Here note that
the core circle of $\Gamma$ is parametrized by  $\varphi_2 \in [0, 2\pi]$ for
$r_2=1$ (cf. Figure~\ref{solidtorus}).

\v
\noindent Then we observe that  $\gamma = \frac{1}{2} ( r_1^2 d \varphi_1 +
r_2^2 d \varphi_2)$ is a primitive of $\omega_0$ on $\mathbb{R}^4$
and that $${\gamma |}_{\Gamma} = \frac{1}{2} ( g^2(r_2) d \varphi_1 +
r_2^2 d \varphi_2)= \frac{r_2^2}{2} ( \frac{g^2(r_2)}{r_2^2}  d
\varphi_1 +  d \varphi_2)$$ is a contact 1-form on $\Gamma$. This
can be verified by a direct calculation where $r_2 d \varphi_1 \wedge dr_2 \wedge
d\varphi_2$ is a volume form on $\Gamma$. Also observe that the
core circle of $\Gamma$ is transverse to the contact structure
$\ker ({\gamma |}_{\Gamma}) $ since $\gamma(\frac{\partial}{\partial
\varphi_2})|_{\{r_2=1\}} = 0.5$.

\v

\begin{figure}[ht]
  \begin{center}
     \includegraphics{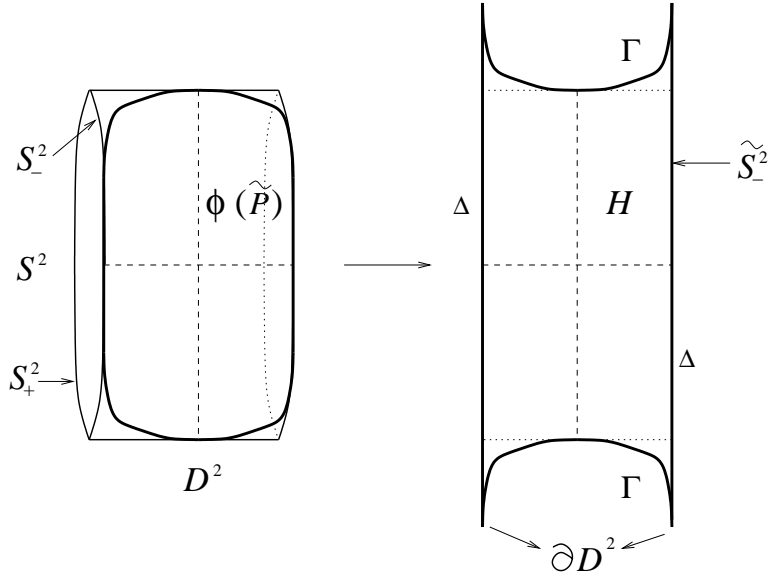}
   \caption{The 2-handle $H$ }
   \label{handle}
    \end{center}
  \end{figure}

\noindent Moreover $(\widetilde{P},
\omega_0)$ is a convex filling of $(\Gamma, {\gamma |}_{\Gamma})$
but we would like to have a concave filling. So we apply the
following trick. We embed $\widetilde{P}$ into a symplectic $S^2
\times D^2$ by a symplectomorphism and take the complement of the
image in $S^2 \times D^2$. Let $(S^2, \sigma_1)$ be a symplectic sphere with area $2
\pi$ and $(D^2, \sigma_2)$ be symplectic disk with area $\pi a^2$.
Denote by $S_{\pm}^2$ the upper and lower hemispheres of area
$\pi$, respectively.  Then $\sigma_1 \oplus \sigma_2$ induces a symplectic form on
$S_{+}^2  \times D^2$. Let $$\phi: P \cong D^2 \times D^2 \to
S_{+}^2  \times D^2 \subset S^2 \times D^2 $$ be a
symplectomorphism. From now on we will identify the symplectic
form on $P$ induced from $\omega_0$ on $\mathbb{R}^4$ with the symplectic form
$\sigma_1 \oplus \sigma_2$ on $S_{+}^2  \times D^2$ by the above symplectomorphism $\phi$.
Define the 2-handle $H$ (see Figure~\ref{handle}) as
$$ H= \overline{S^2 \times D^2 -\phi(\widetilde{P})} .$$

\noindent Now consider the boundary $\partial H$ of the 2-handle $H$. Let us denote
$\phi (\Gamma)$ also by $\Gamma$ by an abuse of notation.
Let $$\Delta = \overline{\partial H \setminus \Gamma} .$$
Observe that $\Delta$ is
fibered by discs $D_x = \widetilde{S_-^2} \times \{x\}$ for $x \in
\partial D^2$, where we have $ S_-^2 \subset \widetilde{S_-^2}
\subset S^2$. This is illustrated in Figure~\ref{handle}: Imagine
the complement of $\phi(\widetilde{P})$ in $S^2 \times D^2$ restricted
to $\partial D^2$. Notice here that $\widetilde{S_-^2}$
is symplectic (with respect to $\omega_0$)  with
fixed area $\frac{7\pi}{4}$ (not $\frac{9\pi}{4}$ as mistakenly
typed in \cite{el-04}) for each $x \in \partial D^2$. This is
precisely because of our identification of $\omega_0$ with
$\sigma_1 \oplus \sigma_2$.

\v

\noindent Next we would like to find an appropriate way to attach this
2-handle $H$ to $Y\times I$ by identifying $\Gamma$ with a
neighborhood $U$ of the binding $B$ of the compatible open
book $\OB_\xi$ in $Y\times \{1 \} \subset Y \times I $. By Giroux
\cite{gi-02},
we can find coordinates $(r, \varphi, u)$ near the binding $B$ of
$\OB_\xi$ such that $$U \cong [0,R] \times (\mathbb{R}^2 / 2\pi
\mathbb{Z}) \times (\mathbb{R}^2 / 2\pi \mathbb{Z})$$ satisfying the
following conditions:
\vskip.1in
(1) $\alpha|_U = h(r) (du + r^2 d\varphi)$ for some positive
function $h$ defined on $[0,R]$ such that $h(r)-h(0)=-r^2$ near
$r=0$, and $h'(r) < 0$ for all $r >0$,
\vskip.1in
(2) $d\alpha$  is symplectic on the pages of $\OB_\xi$ and
\vskip.1in
(3) Pages of $\OB_\xi$ in $U$ are given by $\varphi$=constant.

\v

\begin{figure}[ht]
  \begin{center}
     \includegraphics{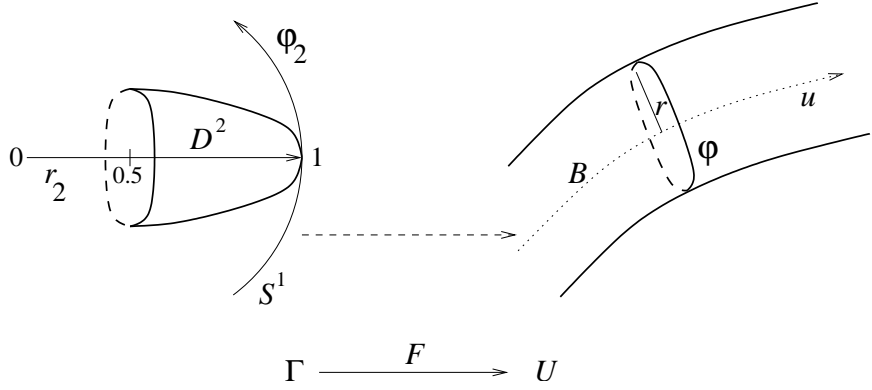}
   \caption{The diffeomorphism $F : \Gamma \to U$ }
   \label{diff}
    \end{center}
  \end{figure}

\noindent Let $a= \frac{R}{2}$ and consider the following map $F : \Gamma
\to U$ (cf. Figure~\ref{diff}) given by the following identifications of coordinates: $$ r=
\frac{g(r_2)}{r_2}, \;\; \varphi=\varphi_1, \;\;u=\varphi_2 . $$

\noindent Notice that
under this map the core circle of $\Gamma$ is mapped onto the
binding $B$. It is clear that $F$ is a diffeomorphism but we would
like $F$ to be a contactomorphism which takes the contact structure
$\ker (\gamma|_\Gamma)$ onto the contact structure $\xi\vert_U$. Well, we will simply choose
our function $g$ (depicted in Figure~\ref{g}) accordingly near $t=1$
so that $F$ becomes a contactomorphism.
The function $$t \to
\frac{g(t)}{t}$$
 is a decreasing function from $[0.5  , 1]$
to $[0, 2a] = [0, R]$. Let $\psi : [0,R] \to [0.5 , 1]$ be
the inverse of this function. Recall that
$$\gamma|_\Gamma = \frac{r_2^2}{2} \big( \frac{g^2(r_2)}{r_2^2} d
\varphi_1 + d \varphi_2 \big).$$

\noindent Hence by the change of coordinates which describes
the diffeomorphism $F$ we get

$$F_* (\gamma|_\Gamma)= (F^{-1})^* (\gamma|_\Gamma)
= \frac{\psi^2 (r)}{2} (r^2 d \varphi +du)$$

\noindent which is a 1-form defined on $U$. Then we have
$ h(r) = \frac{1}{2} \psi^2 (r)$
where $h(r)-h(0)=-r^2$ near $r=0$ so that $ \frac{1}{2} (\psi^2 (r)
- \psi^2 (0)) =-r^2$
and thus $ \psi(r) = \sqrt{1-2r^2}$. Recall that $r_2 = \psi(r)$ and $r=
\frac{g(r_2)}{r_2}$ under the diffeomorphism $F$. Considering that
$\psi : [0,R] \to [0.5, 1]$ is the inverse of the
function $t \to \frac{g(t)}{t}$ we finally obtain
$$r_2 = \big(1-2\frac{g^2(r_2)}{r_2^2}\big)^{1/2}$$

\noindent which implies that $$ g(t) =
\frac{1}{2} t \sqrt{1-t^2}$$

\noindent near $t=1$.  Notice that $g(t)$
has a vertical tangent at $t=1$. This
calculation determines the particular form of the function $g(t)$
near $t=1$. (We calculated the function $g$ slightly
different from the one given in \cite{el-04}.)

\v

\begin{figure}[ht]
  \begin{center}
     \includegraphics{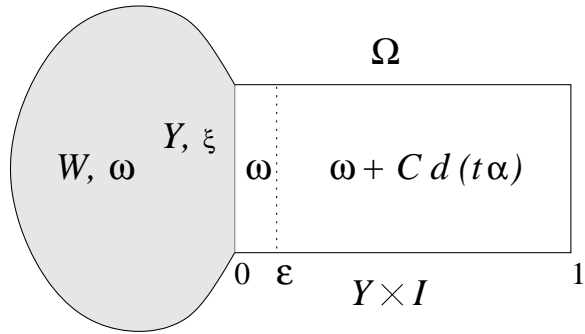}
   \caption{Extending $\omega$ to $\Omega \in H^2(Y\times I)$}
   \label{collar}
    \end{center}
  \end{figure}

\noindent While preparing our 2-handle for gluing we also have to equip the
cobordism $Y \times I$ by an appropriate symplectic form.
It is easy to see (cf. \cite{el-04})
that for an arbitrary constant $C$ there is a symplectic form
$\Omega$ on $Y\times I$ (see Figure~\ref{collar}) which
``extends" $\omega$ and
agrees with $ \omega + C d(t\alpha)$
for $t \in [\varepsilon, 1]$.

\v

\noindent Now observe that $Y \times \{1 \} \subset
Y \times I$ is convex with respect to $d(t\alpha)$. Recall that $(H, \omega_0)$
is a strong concave filling of $(\Gamma, \gamma|_\Gamma)$. Thus by Lemma~\ref{gluing},
we can glue $(H, \omega_0)$ to $(Y \times I , d(t\alpha))$ identifying
$\Gamma$ and $U$ by the contactomorphism $F$ which extends to a symplectomorphism
in some neighborhoods
of $\Gamma \subset H$ and $U \subset Y\times I$. Consequently we have $F^* (d
(t \alpha)) =\omega_0$. Notice that $F^* \omega$ is exact in a
regular neighborhood $\nu(\Gamma)$ of $\Gamma$ in $H$ because the second cohomology
group of $\nu(\Gamma)\cong I \times S^1 \times D^2$ is trivial. Since $F^* \omega$ is exact, there is
a 1-form $\theta$ on $\nu(\Gamma)$ such that $F^* \omega =d \theta $.
Take a smooth cut-off function $\sigma$ on $H$ which vanishes outside
of $\nu(\Gamma)$. Then $d(\sigma \theta)$ defines an extension $\widetilde{\omega}$ of
$F^* \omega$ from $\nu(\Gamma)$ to $H$.

\v

\noindent Finally we are ready to define a symplectic form on $H$ that will
allow us to make this 2-handle attachment in the symplectic
category. Let $\Omega_0 =  \widetilde{\omega} + C \omega_0 $ on $H$
for some constant $C$. It is not hard to see that
$\Omega_0$ will be  symplectic for sufficiently large values of $C$ since
$C^2 \omega_0 \wedge \omega_0 >0$ will dominate the other terms in
$\Omega_0 \wedge \Omega_0$ on a compact manifold.
Here notice that
we have a well-defined symplectic form on $(Y \times I) \cup_{U=
F(\Gamma)} H $ since $\Omega = \omega +Cd(t\alpha)$ on $(Y \times I)$
is identified
with $\Omega_0$ on $H$ in the gluing region. This is because
$\widetilde{\omega}$ is an extension of $F^* \omega$ and $F^* (C d
(t \alpha)) = C \omega_0$ so that $F^* \Omega =\Omega_0$.

\v
\noindent On the other hand, by attaching the 2-handle $H$ we perform a Dehn
surgery on the 3-manifold $Y$ to yield a 3-manifold $Y'$
which fibers over the circle. This should be clear
since we take out a neighborhood $U$
from $Y$ and glue in $\widetilde{S_{-}^2 }\times \partial D^2$ to
cap off each page $F$ of $\OB_{\xi}$ by a disk $D_x = \widetilde{S_{-}^2 }\times \{x\}$.
Let $\widehat{F}$ denote the closed surface obtained by capping off
a page $F$ by gluing a 2-disk $D_x$
along its boundary. Consider the 2-form $\Omega\vert_Y = \Omega\vert_{Y \times \{1\}} = \omega+ C d \alpha $.
We know that $d \alpha$ is symplectic on every page $F$ of $\OB_{\xi}$.
Thus $\omega+ C d \alpha$ will be a symplectic form on $F$ for
sufficiently large values of $C$.
\v
\noindent
Recall that we identified the symplectic forms $\Omega$
and $\Omega_0$ when we attached the symplectic 2-handle $H$.
Also note that $D_x$ is symplectic with respect
to $\omega_0$. Consequently
since every page $F$ of $\OB_{\xi}$ is symplectic (with respect to $\Omega$)
and the disk $D_x$ is symplectic as well (with respect to $\Omega_0$)
we get a fibration over $S^1$ for which
$\omega' = \Omega_0|_{Y'}$ restricts to a symplectic form on each fiber
$\widehat{F}$ for sufficiently large values of $C$. We will call such a surface
fibration over $S^1$ a \emph{symplectic} fibration over $S^1$.
Note that we have the freedom
to choose $C$ as large as we wish. Also
note that in order to prove that we have
a symplectic fibration over $S^1$ after surgery
we had to use the compatibility of $\xi$ and $\OB_{\xi}$.

\v
\noindent
Denote by $(W', \omega')$
the resulting symplectic 4-manifold obtained by attaching the symplectic
2-handle $H$ to the given weak filling $(W, \omega)$ of $(Y, \xi)$.
To finish Eliashberg's construction we need to cap off the
symplectic fibration $\partial W'=Y' \to S^1$ by a symplectic 4-manifold.
Let $\phi$ be the topological monodromy of this surface fibration. Then we can \emph{smoothly}
fill in $ - Y'$ (see \cite{ao-02}) by a
symplectic Lefschetz  fibration over $D^2$ with regular fiber $\widehat{F}$
since the monodromy $\phi^{-1}$ of $-W'$ can be written as a product of
\emph{right-handed} Dehn twists
by Lemma~\ref{positive}.

{\Lem [\cite{ao-02}] \label{positive} Any element in Map$(\widehat{F})$ can be
expressed as a product of non-separating right-handed Dehn twists.}

\begin{proof}
Recall that the relation (cf. Lemma~\ref{rel})
$$
{(t_{a_1} t_{a_2} \cdots t_{a_{2g}})}^{4g+2} =t_c
$$
in $\Gamma_F$ induces a relation
$$
{(t_{a_1} t_{a_2} \cdots t_{a_{2g}})}^{4g+2} =1
$$
in $\Gamma_{\widehat{F}}$. We conclude that
$t_{a_1}^{-1}$ is a product of non-separating right-handed Dehn twists.
Therefore any left-handed non-separating Dehn twist --- being conjugate to
$t_{a_1}^{-1}$ --- is a product of non-separating right-handed Dehn twists.
This finishes the proof of the lemma combined with the fact
that $\Gamma_{\widehat{F}}$ is generated by (right and left-handed) non-separating Dehn twists.
\end{proof}

\v
\noindent
In fact, Eliashberg proves a ``symplectic"
version of Lemma~\ref{positive} in \cite{el-04}
so that we can actually fill in $ -\partial W'= -Y'$ \emph{symplectically}
by a symplectic 4-manifold. The point here is that when we measure the topological
monodromy of a symplectic fibration  $Y' \to S^1$
we do not take into account the symplectic
structure on the fiber. But to fill in such a symplectic fibration
\emph{symplectically} we need to
measure the holonomy (i.e., ``symplectic"  monodromy) of this fibration, which we describe below.
Suppose that the symplectic fibration
$Y' \to S^1$ is normalized so that $\int_{\widehat{F}} \omega'=1.$
Since the 2-form $\omega'$ is positive on the fibers its kernel $\ker \omega'$ is
a 1-dimensional line field on $Y'$ transverse to the fibers.
The flow generated by a vector field which directs this line field determines a holonomy
automorphism $\mbox{Hol}(\omega') \colon \widehat{F}_0 \to \widehat{F}_0$ of
a fixed fiber $\widehat{F}_0$. This is an area and orientation preserving
diffeomorphism (i.e., a symplectomorphism)  which defines
$(Y', \omega')$ uniquely up to fiber preserving diffeomorphism fixed
on $\widehat{F}_0$.
\v
\noindent
Now let $(V, \eta)$ denote the symplectic Lefschetz  fibration over $D^2$ mentioned above
with regular fiber $\widehat{F}$ which will be used to fill in the symplectic
fibration $-Y' \to S^1$. Since $(V, \eta) \to D^2$ is a symplectic Lefschetz fibration,
the symplectic 2-form $\eta$ restricts to a symplectic form on each regular fiber and moreover
we can assume that $\eta\vert_{\partial V}$ integrates
to $1$ on the fibers of the symplectic
fibration $\partial V \to S^1$.
If we can choose $(V, \eta)$ such that $\mbox{Hol}(\omega'\vert_{Y'})^{-1}= \mbox{Hol}(\eta\vert_{\partial V})$
then we are done since we can glue $(W', \omega')$ to $(V, \eta)$ symplectically.
Eliashberg constructs such a symplectic Lefschetz fibration over $D^2$ in \cite{el-04}.
In fact as it is shown in \cite{km-04}
that it suffices to prove Lemma~\ref{holo} below.
(See also Section~\ref{last} for another argument for
the sufficiency of Lemma~\ref{holo}.)


{\Lem [Kronheimer-Mrowka \cite{km-04}] \label{holo}
Let $\Sigma$ be a closed symplectic surface of area $1$ and genus $g>1$.
Let $\phi : \Sigma \to \Sigma$ be an area preserving diffeomorphism that
is smoothly isotopic
to the identity. Then there is a symplectic Lefschetz fibration
$p : (V, \eta) \to D^2$ such that $p^{-1} (1) = \Sigma$ and
$Hol(\eta \vert_{\partial V})= \phi$.}

\v
\noindent As it is pointed out in \cite{el-04}, we could alternatively
use Lemma~\ref{surfacebund} to cap off a symplectic fibration
over $S^1$ by a symplectic surface bundle over a surface with boundary.
Recall that a group $G$ is said to be perfect
if it is equal to its commutator subgroup $[G,G]$. In other words, $G$ is perfect if
and only if every element in $G$ can be expressed as a product of commutators.
Yet another way of characterizing the perfectness of a group is
given by the triviality of its first homology group $H_1(G) = G/[G,G]$.

\v
\noindent It is well-known that the
mapping class group of a surface of genus greater than two is perfect.
This is a consequence of the lantern relation (cf. \cite{jo-79}) in the mapping class groups
which essentially says ``three equals four".
The fact that one can smoothly fill
in a smooth surface bundle over $S^1$ by a smooth surface (of genus $> 2$) bundle over a surface
with boundary easily follows from the perfectness of the corresponding mapping class group (of genus $> 2$).
Here we need a symplectic version of this fact which is provided by Kotschick
and Morita \cite{komo-05}. Let
$Symp_\sigma \Sigma$ denote the  group of all symplectomorphisms of the closed
symplectic surface $(\Sigma, \sigma)$ with respect to a prescribed symplectic
form $\sigma$ on $\Sigma$ which is normalized such that $\int_{\Sigma}\sigma =1.$

{\Lem [Kotschick-Morita \cite{komo-05}] \label{surfacebund} If the genus of $\Sigma$
is greater than two then $Symp_\sigma \Sigma$ is perfect. }

\v
\noindent The restriction in Lemma~\ref{surfacebund}
on the genus of the fiber is not a serious one
since in the construction above one can arbitrarily increase the genus of
the page of $\OB_{\xi}$ (which is compatible with $(Y, \xi)$) by
positively stabilizing $\OB_{\xi}$ (cf. \cite{gi-02}) to begin with.

\subsection{Etnyre's construction:} In Section~\ref{strongembed} we showed that
to find an embedding of a strong filling one can use an embedding of a
Stein filling. Etnyre's idea in \cite{et-04} was to find an
embedding of a weak filling using an embedding of
a strong filling.
Suppose that $(W, \omega )$ is a weak filling of a contact 3-manifold $(Y, \xi)$.
Etnyre showed that $(W, \omega )$ can be embedded into
a symplectic 4-manifold $(W', \omega ')$ which weakly fills its
boundary $(\partial W'=Y', \xi ')$, where $Y'$ happens to be a
integral homology sphere. Now by a homological argument the
symplectic structure $\omega '$ can be perturbed near the boundary
so that $(W', \omega ')$ strongly fills $(Y', \xi ')$. Therefore $(W',
\omega ')$ can be embedded into a closed symplectic 4-manifold
$(X, \omega_X)$ by Proposition~\ref{stembed} and hence $(W,
\omega ) \subset (W', \omega ') $ can be embedded symplectically
into $(X, \omega_X)$. Below we proceed with the details.

\v
\noindent
Let $(W, \omega )$ be a weak filling of $(Y, \xi )$
and let $\OB_\xi$ be an open book compatible with $(Y, \xi )$. We
can assume that the binding $B$ of $\OB_\xi$ is connected. Let
$\phi$ be the monodromy of this open book and use Lemma~\ref{factor} to express
$\phi$ as $$ \phi = t_c^m \circ  t_{\gamma_1}^{-1}\circ \cdots \circ
t_{\gamma_n}^{-1} . $$ Now Legendrian Realize $\gamma_n$ (cf. Lemma~\ref{lrp})
on a convex page of $\OB_\xi$ and perform contact $(-1)$-surgery on $\gamma_n$. The new
open book will have monodromy $$ \phi\circ t_{\gamma_n} = t_c^m \circ
t_{\gamma_1}^{-1} \circ \cdots \circ t_{\gamma_{n-1}}^{-1} .$$ Repeat this for
all the curves $\gamma_i$ (for $i=n-1, \cdots, 1$) to get down to $t_c^m$ as
in the proof of Proposition~\ref{conca}.
Denote by $(Y', \xi')$ the contact 3-manifold obtained as a
result of the contact $(-1)$-surgeries above. Then  $(Y', \xi')$ is
compatible with the open book whose monodromy is given by $t_c^m$, by Lemma~\ref{compa}.

\v

\noindent Recall that by Theorem~\ref{strong} we can attach a symplectic
2-handle to a strong filling along a Legendrian knot in its convex boundary in such
a way that the symplectic structure extends to the
2-handle and the new symplectic $4$-manifold strongly fills its boundary. In this gluing
process, however, the Liouville (i.e., symplectically dilating)  vector field
is used only in a neighborhood
of the attaching circle. It turns out that if $L\subset (Y, \xi )$ is Legendrian and $(W,
\omega )$ is a weak filling of $(Y, \xi )$ then there is always a symplectically dilating
vector field near $L$, implying

{\Prop [Ding-Geiges \cite{dg-01}] \label{weak}
Suppose that $(Y', \xi ')$ is given by contact $(-1)$-surgery along $L\subset (Y, \xi )$.
If  $(Y, \xi )$ is weakly fillable then so is $(Y', \xi ')$. }

\v
\noindent
Hence there exists a weak filling $(W', \omega')$ of $(Y',
\xi')$ obtained by attaching symplectic 2-handles to $(W, \omega)$.
The page of the compatible open book $\OB_{\xi'}$ is a genus $g$ surface with one
boundary component. Consider the curves $a_i$ depicted in
Figure~\ref{surface}. Legendrian realize $a_i$'s and perform
contact $(-1)$-surgery on each $a_i$ to get $(Y'', \xi'')$ compatible with the
open book $\OB_{\xi''}$ whose monodromy is given by $$t_c^m \circ  t_{a_1}^{-1}\circ  \cdots \circ
t_{a_{2g}}^{-1} .$$ It is not hard to see that $Y''$ is an integral
homology sphere. Moreover, by Proposition~\ref{weak}, there exists a weak filling $(W'', \omega'')$
of $(Y'', \xi'')$. Then we use Proposition~\ref{modif} to modify the symplectic form
$\omega''$ near the boundary so that it is a strong filling of
$(Y'', \xi'')$. Note that  $(Y'', \xi'') $ has a concave filling
by Proposition~\ref{conca}. Thus we cap
off $(W'', \omega'')$ by this concave filling using Lemma~\ref{gluing}
to get a closed
symplectic 4-manifold $(X, \omega_X)$ in which $(W, \omega)$ sits as a symplectic subdomain.

\v
{\Prop  [Eliashberg \cite{el-91}, \cite{el-04}; Ohta-Ono \cite{ohon-99}]\label{modif}
Any weak filling of a rational homology sphere can be
deformed into a strong filling by modifying the symplectic form
near the boundary.}

\v
\noindent
The main step in Etnyre's construction is embedding a
weak filling of an arbitrary contact 3-manifold into
a weak filling of an integral homology sphere.
We would like to point out here that this follows also from a result that
was obtained by Stipsicz in \cite{st-03}.
Namely, Stipsicz showed the existence of a Stein cobordism from an arbitrary
contact 3-manifold to an integral homology sphere. Stipsicz's
construction (which we describe below) can be slightly modified to imply the main step above.

\begin{figure}[ht]
  \begin{center}
     \includegraphics{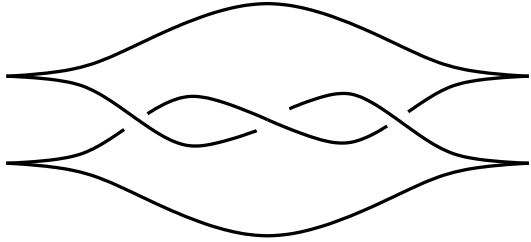}
   \caption{A right-handed Legendrian trefoil knot $K$}
   \label{tref}
    \end{center}
  \end{figure}

\v
\noindent
Let $(W, \omega )$ be a weak filling of $(Y, \xi )$.
Consider the right-handed Legendrian trefoil knot $K$ as depicted in Figure~\ref{tref}
in the standard contact $S^3$,
having $\tb (K)=1$. To construct such a cobordism start with a contact surgery diagram $\Li $ of $(Y,
\xi )$ and for every knot $L_i$ in $\Li$ add a copy $K_i$ of $K$ into the diagram
linking $L_i$ once, not linking the other knots in $\Li$. Adding
symplectic 2-handles along $K_i$ we get $(W', \omega')$ and the resulting
3-manifold $Y'$ is an integral homology sphere. To see this just convert the
contact surgery diagram into a smooth handlebody diagram and calculate the first homology.
Observe that the topological framing of $K$ is $0$. Denote by
$\mu_i$ a small circle meridional to $K_i$ and ${\mu}^{\prime}_i$ a small circle meridional
to $L_i$ for $i =1, \cdots, n$. Recall that
$H_1(Y', \mathbb{Z})$ is generated $[\mu_i]$ and $[{\mu}^{\prime}_i]$ and the relations are
$[{\mu}^{\prime}_i]=0$ and $[\mu_i] + \sum_{j \neq i} lk (L_i, L_j) [{\mu}^{\prime}_j] =0$.
It follows that
$H_1(Y', \mathbb{Z})=0$.

\v
\noindent
Although it was not considered in \cite{st-03}, Stipsicz's construction
immediately implies the main step above because one can add a symplectic 2-handle along a
Legendrian knot in the boundary of a weak filling to extend it to another
weak filling by Proposition~\ref{weak}.

\section{A hybrid solution}\label{last}
\v
\noindent In this section we suggest another symplectic embedding of
a weak filling into a closed symplectic 4-manifold which is obtained by
a mixture of the ideas we discussed so far. First we note
that it is possible to attach symplectic 1-handles (as well as symplectic 2-handles)
to a weak filling to extend it to another weak filling.
Suppose that $(W, \omega )$ is a weak filling of a contact 3-manifold $(Y, \xi)$.
Now we proceed
as in the second proof of Proposition~\ref{conca} to embed $(W, \omega)$ into
a weak filling $(W', \omega')$ by attaching symplectic 1- and 2-handles so that
the resulting contact structure on the boundary $\partial W'$ has a compatible
open book whose page has only one boundary component and whose monodromy
is just one right-handed boundary-parallel Dehn twist. Then we attach a symplectic
2-handle to $(W', \omega')$ along the binding of this open book and we get a
symplectic fibration over a circle with topologically trivial monodromy
on the other end of the cobordism given by this
2-handle attachment. Finally we cap off this surface bundle by a
symplectic Lefschetz fibration over $D^2$ using Lemma~\ref{holo}.

\section{Final comments}\label{hybrid}
\v
\noindent
The presentation in this article may suggest that
Eliashberg's method is unnecessarily long but he constructs from scratch
a symplectic 2-handle
to be attached to the binding of a compatible open book---which is original.
It would be very interesting to
interpret this \emph{new} symplectic surgery
in terms of contact surgery. Unfortunately,
there does not seem to exist a natural contact
structure on the symplectic fibration over $S^1$ obtained by this surgery.
This is exactly the point where Eliashberg's method differs from
the method of Etnyre. In Etnyre's construction one always makes use of
the contact structures on the boundaries of symplectic 4-manifolds
to glue them symplectically. In fact, based on Giroux's correspondence,
Etnyre mostly deals with open books compatible with these
contact structures rather than the contact structures directly.
In Eliashberg's construction, however, at one point or another
we have to glue a
symplectic (Lefschetz) fibration to a symplectic 4-manifold whose boundary
symplectically fibers over $S^1$. This is achieved by matching up the
holonomy diffeomorphisms on the boundaries and contact structures
are not visible
in this picture. It might be
worth pointing out that the proof of the
non-triviality of the contact Heegaard Floer invariant
of a fillable contact structure
follows from Eliashberg's embedding but it is not clear whether or not
it follows from Etnyre's construction.

\v
\noindent
Also it is intriguing to note that most
of the discussion above relies on the relation
$$t_c = ( t_{a_1} \circ  t_{a_2}
\circ \cdots \circ  t_{a_{2g-1}} \circ t_{a_{2g} })^{4g+2} $$
given in Lemma~\ref{rel} in the mapping class group of a
surface $F$ with one boundary component which implies
$$1 = ( t_{a_1} \circ  t_{a_2}
\circ \cdots \circ  t_{a_{2g-1}} \circ t_{a_{2g} })^{4g+2} $$
in the mapping class group of $\widehat{F}$, where $\widehat{F}$ denotes
the closed surface obtained by
capping off the surface $F$ by gluing a
2-disk along $\partial F$. This latter relation says that identity can be
expressed as a product of right-handed Dehn twists in the mapping class group of a closed surface.
It is not possible, however,
to express the identity as a product of right-handed Dehn
twists in the mapping class group of a surface with non-empty
boundary (cf. \cite{os-04}).

\v
\noindent {\bf {Acknowledgement}}: We would like to thank Selman Akbulut for encouragement,
Tolga Etg\"{u} for helpful
conversations and Andr\'as Stipsicz for reading a preliminary version and
suggesting numerous improvements. The author was partially supported by the Turkish
Academy of Sciences.

\end{document}